\documentclass[a4paper,11pt]{article}
\pdfoutput=1 

\usepackage{jheppub} 

\usepackage{comment}

\numberwithin{equation}{section}
\gdef\ffrac#1#2{\textstyle\frac{#1}{#2}\displaystyle}
\gdef\be{\begin{equation}}
\gdef\ee{\end{equation}}
\gdef\e{\epsilon}
\gdef\a{\alpha}
\gdef\p{\partial}

\gdef\d{{\delta}}
\gdef\l{{\lambda}}
\gdef\i{{\mathrm i}}
\gdef\erm{{\mathrm e}}

\gdef\b{{\beta}}

\title{\boldmath $T\overline T$-deformed modular forms}

\author[a,b]{John Cardy}

\affiliation[a]{Department of Physics, University of California, Berkeley CA 94720, USA}
\affiliation[b]{All Souls College, Oxford OX1 4AL, UK}

\emailAdd{cardy@berkeley.edu}

\abstract
{Certain objects of conformal field theory, for example partition functions on the rectangle and the torus, and one-point functions on the torus, are either invariant or transform simply under the modular group, properties which should be preserved under the  $T\overline T$ deformation. The formulation and proof of this statement in fact extents to more general functions such as $T\overline T$ deformed modular and Jacobi forms. We show that the deformation acts simply on their Mellin transform, multiplying it by a universal entire function. Finally we show that Maass forms on the torus are eigenfunctions of the $T\overline T$ deformation.}

\begin{document} 
\maketitle

\section{Introduction and statement of results}\label{s1}
\subsubsection*{Holomorphic modular forms.}
Modular forms and their cousins play an essential role both in mathematics, for example elliptic curves and number theory, and in mathematical physics, \em e.g. \em conformal field theory (CFT) and integrable lattice models. In general, they are functions, defined in the upper half $\tau$-plane, of the form
\be\label{e1.1}
F_1(\tau)=\sum_{n=0}^\infty a_nq^{\Delta+n}\\,
\ee
where $q\equiv \erm^{2\pi\i\tau}$, $a_0\not=0$, and the sum converges in $|q|<1$. Moreover they have simple transformation rules under the generators $S:\tau\to-1/\tau$ and $T:\tau\to\tau+1$ of the modular group $\Gamma_1=$SL$(2,{\mathbb Z})$:
\be
F_1(-1/\tau)=(-\i\tau)^kF_1(\tau)\,,\qquad F_1(\tau+1)=\erm^{-2\pi\i\Delta}F_1(\tau)\,.
\ee
For $\Delta$ a positive integer or zero, and $k$ an even integer, $F_1$ is a modular form of weight $k$, while for $k=0$ and $\Delta$ a negative integer it is a modular function, invariant under $\Gamma_1$, but here we do not make such restrictions, allowing, for example, arbitrary powers of the Dedekind function 
$\eta(\tau)=q^{1/24}\prod_{n=1}^\infty(1-q^n)$. We also include Jacobi forms, which depend on a second variable and which transform like theta functions. 
 
More generally we may consider a vector space of such functions transforming according to some representation of $\Gamma_1$ (or a subgroup, not considered here), examples being the characters of a chiral algebra in a CFT.

In what follows, however it is more useful to regard these as functions of $\d\equiv-\i\tau$, with $q=\erm^{-2\pi\d}$ and $\d\in{\mathbb H}\equiv\{{\rm Re}\,\d>0\}$, so that $F(1/\d)=\d^kF(\d)$.

For all such functions, we shall show that there is a family of deformations labeled by a real parameter $\a>0$ which satisfies

\noindent{\bf Theorem 1 (deformed holomorphic forms.)} \em Given a form $F_1(\d)=\sum_{n=0}^\infty a_n\erm^{-2\pi(\Delta+n)\d}$
with the property that $F_1(1/\d)=\d^kF_1(\d)$, the deformed function
\begin{equation}\label{e1}
F_1^\a(\d)=\sum_{n=0}^\infty a_n\frac{(1+\sqrt{1+8\pi(\Delta+n)\a\d})^{1-k}}{\sqrt{1+8\pi(\Delta+n)\a\d}}\,\erm^{-(1/2\a)(\sqrt{1+8\pi (\Delta+n)\a\d}-1)}
\end{equation}
also satisfies $F_1^{\a}(1/\d)=\d^k\,F_1^{\a}(\d)$,
as long as both sides converge, which, for $\Delta<0$, restricts $8\pi|\Delta|\a<{\rm Re}\,\d<(8\pi|\Delta|\a)^{-1}$.\em

Note that $\lim_{\a\to0}F_1^\a=2^{1-k}F_1$.
It is instructive to rewrite (\ref{e1}) in terms of $\beta=2\pi\a\d$ and $\d$, so that the exponential is $q^{(\Delta+n)_\beta}$, where we have introduced the $\beta$-deformed
number
\be\label{xb}
x_\beta\equiv(1/2\beta)(\sqrt{1+4\beta x}-1)\,.
\ee
Note that 
\be
(\beta x_\beta)^2+\beta x_\beta=\beta x\,,
\ee
so that $\beta x_\beta$ transforms according to a nonlinear representation of the additive semigroup on $\{\beta\in\overline{\mathbb R}_+\}$, isomorphic to the semigroup $\mathfrak T$ generated by the $T\overline T$ flow.

\subsubsection*{Dirichlet series and Mellin transform.}
To every form (\ref{e1.1}) with $\Delta>0$ may be associated a Dirichlet series
\be
\phi(s)=\sum_{n=0}^\infty\frac{a_n}{(\Delta+n)^s}\,,
\ee
which converges for Re\,$s>k$ and can be continued outside this region to obey a simple reflection relation under $s\to k-s$. This also simply related to the Mellin transform
$R(s)\equiv\int_0^\infty\d^{s-1}F(\d)d\d$ by $\phi(s)=((2\pi)^s/\Gamma(s))R(s)$. The modular property  of $F$ then implies that the integral converges for Re\,$s>k$ with a continuation having the reflection property $R(s)=R(k-s)$. 

One may then ask how the series associated to the deformed form is related to that associated to the seed form (\ref{e1.1}). However, after the deformation, $F^\a$ as given in (\ref{e1}) no longer has the form of a series in powers of $q$ when expressed in terms of $\a$ and $q$, but rather in terms of $\beta=\a\d$ and $q$. Thus the deformed Dirichlet series defined by the substitution $\Delta+n\to(\Delta+n)_\beta$ is no longer simply proportional to the Mellin transform of $F^\a$.

 It turns out that it is the latter which enjoys simple properties, given by

\noindent{\bf Theorem 1a (deformed Mellin transform).} \em The Mellin transform $R^\a(s)$ of the $T\overline T$ deformation of a modular form of degree $k$ is related to that of the undeformed form by\em
\be
R^\a(s)= I^\a(k,s)\,R^0(s)\\,
\ee
\em where $I^\a(k,s)$ is a universal entire function of $s$, satisfying $I^\a(k,s)=I^\a(k,k-s)$ so that $R^\a(s)$ inherits the reflexion property and the zeroes of $R^0(s)$.\em

Thus the Mellin transform effectively diagonalizes the $T\overline T$ deformation.

\subsubsection*{Real analytic forms.}
Apart from the above singularities when  $\Delta<0$, $F_1^\a(\d)$ is holomorphic in Re\,$\d>0$, but clearly the periodicity under $T:\d\to\d-\i$ is lost due to the irrationality of the exponents. However, at the cost of losing holomorphicity, there is a generalization to the full modular group: 

\noindent{\bf Theorem 2 (deformed real analytic forms.)} \em Given a real form on the half-plane
\be\label{1.7}
F_2(\d)=\sum_{n=0}^\infty\sum_{p\in{\mathbb Z}} a_{n,p} \erm^{-2\pi(\Delta+n)\d_1+2\pi\i p\d_2}
\ee
with $a_{n,p}=a^*_{n,-p}$, which satisfies
$F_2(1/\d)=|\d|^kF_2(\d)$ and $F_2(\d+\i)=F_2(\d)$,
then
\begin{eqnarray}\label{e2}
F_2^\a(\d)=
\sum_{n=0}^\infty\sum_{p\in{\mathbb Z}} a_{n,p}&&\frac{(1+\sqrt{1+8\pi(\Delta+n)\a\d_1+(4\pi p\a\d_1)^2})^{1-k}}{\sqrt{1+8\pi(\Delta+n)\a\d_1+(4\pi p\a\d_1)^2}}\nonumber\\
&&\times \erm^{-(1/2\a)(\sqrt{1+8\pi (\Delta+n)\a\d_1+(4\pi p\a\d_1)^2}-1)+2\pi\i p\d_2}\,,
\end{eqnarray}
satisfies
\be\label{eth2}
F_2^\a(1/\d)=|\d|^kF^\a_2(\d)\quad\mbox{and}\quad F^\a_2(\d+\i)=F^\a_2(\d)\,.
\ee\em

In the case $k=0$, when $F_2$ is modular invariant, we have the alternative and simpler

\noindent{\bf Theorem 2a.} \em Given a real form as in Thm.~2 but satisfying  
$F_2(1/\d)=F_2(\d)=F_2(\d+\i)$, that is a modular invariant,
then
\be\label{e2a}
F_{2a}^\a(\d)=
\sum_{n=0}^\infty\sum_{p\in{\mathbb Z}} a_{n,p}\erm^{-(1/2\a)(\sqrt{1+8\pi (\Delta+n)\a\d_1+(4\pi p\a\d_1)^2}-1)+2\pi\i p\d_2}
\ee
is also modular invariant.\em

Again, Thms.~2 and 2a hold for $\a>0$ and, if $\Delta<0$, $8\pi|\Delta|\a<{\rm Re}\,\d<(8\pi|\Delta|\a)^{-1}$.

\subsubsection*{Maass forms.}
Another genre of real-valued functions over $\mathbb H$ are Maass (cusp) forms, whose  main characteristic is that, in addition to being $\Gamma_1$ invariant, they are eigenfunctions of the invariant Laplacian on the fundamental region ${\cal F}=\mathbb H/\Gamma_1$. We shall show that there is a close relation between the $T\overline T$ deformation and Maass forms, in fact we have

\noindent{\bf Theorem 3.} \em Any Maass form is invariant (up to a multiplicative constant) under the $T\overline T$ semigroup $\mathfrak T$. \em

\subsubsection*{Origins of this work.}
Although we shall rigorously establish these results, they are motivated by non-rigorous arguments based on examples drawn from recent work in theoretical physics on the so-called $T\overline T$ deformation of a two-dimensional CFT. This  began with the paper of Zamolodchikov\cite{Zam1}, and since then differing but related explanations have been put forward for its solvability: in terms of a coupling to random geometry\cite{JC1}, as a particular form of quantum gravity\cite{DubJT}, through its holographic interpretation, or as a state-dependent diffeomorphism\cite{Conti1}.  Since it has an explicitly locally rotationally invariant form (in euclidean signature) it is expected to preserve the modular properties of the original theory defined in domains such as the torus, or, more simply, a rectangle. However, the formalisms mentioned above tend to obscure this symmetry, and it was in formulating a proof that it does in fact hold that the author realized that the arguments apply to more general mathematical objects beyond those which arise in CFTs, such as (\ref{e1.1}). 
In a companion physics paper we shall amplify this interpretation and show that the resulting equations are those of a non-interacting fluid.)
\subsubsection*{Other related work.}
The deformations of modular forms discussed here are quite different from those arising in string theory amplitudes in a gravitational plane wave background, first computed in \cite{pw} and studied in generality in \cite{Berg}. From a world sheet perspective, the latter is a massive deformation, relevant in the infrared, while $T\overline T$ is relevant in the ultraviolet. More explicitly, its effect, for example on the $\eta$-function is to make
\be
\prod_{n=1}^\infty(1-q^n)\to\prod_{n=1}^\infty(1-q^{\sqrt{m^2+n^2}})\,,
\ee
while the effect of $T\overline T$ is to modify the power of each term its series expansion
\be
\prod_{n=1}^\infty(1-q^n)=\sum_{k=0}^\infty\chi_kq^k\to
\sum_{k=0}^\infty\chi^\beta_kq^{(1/\beta)(\sqrt{1+2\beta k}-1)}\,.
\ee

The recent physics literature on $T\overline T$ amounts to several hundred papers. 
Those most relevant  to the present discussion are listed as Refs.~(\cite{Zam2}--\cite{dattach}) 
Recently Benjamin \em et al.\em \cite{benj} extensively treated CFT torus partition functions from the point of view of harmonic analysis, emphasizing the role of Maass forms, without, however, noting the connection to the $T\overline T$ deformation. 
\subsubsection*{Outline.}
The outline of this paper is as follows. In Sec.~\ref{s2} we briefly describe the $T\overline T$ deformation in a non-rigorous manner, re-interpreting it as a coupling of the CFT to an elastic medium, and applying it first to the partition function in a rectangle, then to a 1-point function on the torus. We argue that these are most easily understood through their Laplace transforms, which give a complexified version of the microcanonical ensemble. (In a companion physics paper we shall amplify this interpretation and show that the resulting equations are those of a non-interacting fluid.) These examples are sufficiently general as to give the basis for a rigorous definition of the deformation of modular forms on $\Gamma_1$ in general. However, this section is not necessary for the remainder of the mathematical discussion, which is self-contained. In Secs.~\ref{s3} and \ref{s5} we give proofs of Thms.~(1,2,2a) by showing that the above definitions are equivalent both to the series expansions in the statement of the theorems, and also imply that $F_{1,2}^\a$ are related to their undeformed counterparts by integral transforms which preserve their modular properties. In Sec.~\ref{s4} we discuss Dirichlet series, and in Sec.~\ref{s6} the relation to Maass forms. Finally in Sec.~\ref{s7} we discuss some examples.

\section{Physics motivation}\label{s2}
\subsection{$T\overline T$ as a coupling to an elastic medium}\label{s2.1}

Briefly, the $T\overline T$ deformation of a given CFT in flat space is a family of non-local field theories ${\cal T}^\l$, parametrized by a real number $\l$, in which the infinitesimal flow ${\cal T}^\l\to{\cal T}^{\l+\d\l}$ is formally defined in the path integral representation by adding a perturbation
\be
\d\l\int\det T^\l(x)d^2x=\ffrac12\d\l\e^{ik}\e\e^{jl}\int T^\l_{ij}(x)T^\l_{kl}(x)d^2x
\ee
to the action, or, equivalently, inserting it into correlation functions. Here $T^\l(x)$ is the local energy-momentum, or stress, tensor of the deformed theory, assumed to exist. $x=(x_1,x_2)$ are cartesian coordinates and we use the summation convention. By definition this generates a semigroup action on the space of deformed theories, isomorphic to the additive group on $\overline{\mathbb R}_+$.

The effect of this quadratic perturbation in the path integral may be written as an integral over a symmetric tensor field $\varepsilon$:
\be\label{ela}
\propto\int[d\varepsilon_{ij}]\erm^{\int\varepsilon_{ij}T^\l_{ij}d^2x +(1/2\d\l)
\e^{ik}\e^{jl}\int \varepsilon_{ij}(x)\varepsilon_{kl}(x)d^2x}\,,
\ee
which, being gaussian, is given by the value of the exponent at the saddle
$\varepsilon_{ij}=-(\d\l)\e^{ik}\e^{jl}T^\l_{kl}$. The fact that $T^\l_{kl}$ is conserved and symmetric then implies that $\varepsilon_{ij}$ may be written as $\frac12(\p_iu_j+\p_ju_i)$. 

We may interpret the field $u_j(x)$ as the displacement of a particle, initially at $x$, in an elastic solid, with $\varepsilon_{ij}(x)$ being the strain and the second term in (\ref{ela}) the elastic energy.  

Integrating along a contour $C[X,x]$ from a fixed point $X$  to $x$,
\be\label{CXx}
\p_\l u_i(x)=-\e^{ik}\e^{jl}\int_X^x T^\l_{kl}(x')dx'_j\,,
\ee
where the integral may be recognized as the flux $N^\l_k$ of the stress current $T^\l_{k.}$, that is the total force acting across $C[X,x]$, and is independent of the contour since the current is conserved. Thus the separation $R^\l(a,b)$ between two points initially at $x_a$ and $x_b$ satisfies
\be\label{rn}
\p_\l R_i^\l(a,b)=-\e_{ik}N_k^\l(a,b)\,.
\ee
This interpretation of $u(x)$ as a dynamical field in a fixed frame avoids the paradoxes which may arise in thinking of $x\to x+u(x)$ as a field-dependent diffeomorphism\cite{Conti1} in the quantized theory, and there is no requirement of general covariance.

In general (\ref{rn}) relates two fluctuating quantities, but in some situations $C$ is macroscopic and we may consider a statistical ensemble in which $N_k^\l(a,b)$ is fixed, and moreover a protocol in which it is independent of $\l$. In that case $R_i^\l(a,b)$ evolves linearly with $\l$. We now consider a simple example.

\subsection{Rectangular geometry}\label{s2.2}
Consider a sample of this elastic material coupled to a CFT in the initial shape of an $R_1\times R_2$ rectangle, oriented with its sides parallel to the cartesian axes. Thinking of $x_2$ as imaginary time, the undeformed partition function has a spectral decomposition
\be\label{2.5}
Z^0(R_1, R_2)=\int \erm^{-N_2R_2}\rho^0(R_1, N_2)dN_2\\,
\ee
where the force $N_2$ is the energy in the microcanonical ensemble, and, if we imposed periodic boundary conditions in $x_2$, $R_2$ would be the inverse temperature. Here $\rho^0$ is the density of energy eigenstates, weighted by matrix elements to the initial and final boundary states at $x_2=0$ and $R_2$, and is a sum of delta functions. 

Now take $a$ and $b$ in (\ref{rn}) to be the ends of an interval spanning the rectangle along $x_2=$ constant. This gives the change in the width at this height. 
If $T_{12}=0$, that is there is no shear force at the boundaries, conservation implies that the width change is independent of the height, so the sample remains rectangular.  Thus at fixed $N_2$, $R_1^\l=R_1^0-\l N_2$, so we may write, at least formally\footnote{The change of sign in front of $\l$ is due to $\rho^0$ being a density.}
\be\label{2.6}
Z^\l(R_1, R_2)=\int \erm^{-N_2R_2}\rho^0(R_1+\l N_2, N_2)dN_2\,.
\ee
There are several problems with (\ref{2.6}), one being that $R_1+\l N_2$ may become negative. The other is that although, as a consequence  of the symmetry of the action  under $(x_1,x_2)\to (x_2,x_1)$. we would expect $Z^\l(R_1, R_2)$ to be invariant under
$S:R_1\leftrightarrow R_2$, this is obscured in (\ref{2.6}). 

However, again at least formally, it implies that $Z^\l$ obeys the PDE 
\be
\p_\l Z^\l(R_1, R_2)=-\p_{R_1}\p_{R_2}\,Z^\l(R_1, R_2)
\ee
which is symmetric, although this does not imply that its solution must be even if the initial data are, since the operator on the right hand side is not elliptic.

It is necessary to give mathematical meaning to (\ref{2.5},\ref{2.6}) and the subsequent manipulations. The connection between this example and the more general form in (\ref{e1.1}) is that the undeformed partition function $Z^0$ in a rectangle with boundary conditions $T_{12}=0$ in fact takes the form \cite{kleb}:
\be
Z^0(R_1, R_2)=R_1^{c/4}\eta(\i R_2/R_1)^{-c/2}\\,
\ee
for any CFT of central charge $c$. Here $\eta(\tau)=q^{1/24}\prod_{m=1}^\infty(1-q^m)$ 
with $q=\erm^{2\pi\i\tau}$ is Dedekind's function. This is of the form (\ref{e1.1}) with $k=-c/4$
and $\Delta=-c/48$, but we keep these parameters more general in the following discussion. 

In (\ref{2.5}) $Z^0(R_1,R_2)$ is the Laplace transform of $\rho^0(R_1,N_2)$. 
However, for our purposes it is better to do the reverse, that is define
\be
\omega^0(R_1,s)=\int_0^\infty \erm^{-sR'_2}Z^0(R_1,R'_2)dR'_2\,.
\ee
Since $Z_0\sim \erm^{-2\pi\Delta R'_2/R_1}$ as $R'_2\to\infty$, and $\sim\erm^{-2\pi\Delta R_1/R'_2}$ as $R'_2\to0$, convergence in these limits is is uniform in any closed subset of $R_1\in{\mathbb R}_+$ as long as $\Delta>0$ ($c<0$): $\omega^0$ is then a complex analytic function of $s$ apart from poles along the negative real axis. The case when $\Delta<0$ will be discussed later.

The inverse transform is then
\be
Z^0(R_1,R_2)=\int_{-\i\infty}^{\i\infty}\erm^{sR_2}\omega^0(R_1,s)(ds/2\pi\i)\,,
\ee
and, on pulling back the contour to wrap around the poles, gives $\rho^0(R_1,N_2)$ as twice the imaginary part of $\omega^0(R_1,s)$ at $s=-N_2$.

Comparing with (\ref{2.6}), we therefore have
\be
Z^\l(R_1,R_2)=\int_{-\i\infty}^{\i\infty}\erm^{sR_2}\omega^0(R_1-\l s,s)(ds/2\pi\i)\,,
\ee\,
which still requires a definition of $\omega^0(R_1,s)$ for $R_1\in{\mathbb C}$.
However, for a CFT,
 $Z^0(R_1,R'_2)=R_1^{-k}F_1^0(\d')$ where $\d'=R'_2/R_1$ giving
\be
\omega^0(R_1,s)=\int_0^\infty R_1^{1-k}\erm^{-s\d'R_1}F_1^0(\d')d\d'\,,
\ee
which may be analytically continued in $R_1$.

Thus
\be\label{eq}
Z^\l(R_1,R_2)=\int_{-\i\infty}^{\i\infty}\erm^{sR_2}\int_0^\infty (R_1-\l s)^{1-k}\erm^{-s\d'(R_1-\l s)}F_1^0(\d')d\d'(ds/2\pi\i)\,.
\ee
Finally, writing $Z^\l(R_1,R_2)=R_1^{-k}F_1^\a(\d)$ in terms of dimensionless quantities
$\a=\l/(R_1R_2)$ and $\d=R_2/R_1$, and rescaling  $R_1s\to s$,
\be\label{def}
F_1^\a(\d)=\int_{-\i\infty}^{\i\infty}\erm^{s\d}\int_0^\infty (1-\a\d s)^{1-k}\erm^{-s\d'(1-\a\d s)}F_1^0(\d')d\d'(ds/2\pi\i)\,.
\ee
So far the discussion has lacked rigor. The idea is now to use (\ref{def}) as the \em definition \em of the deformed modular form $F_1^\a$, and to prove both that it yields the expansion in (\ref{e1}) and that it has the same transformation law as $F_1^0$ under $S:\d\to1/\d$. 

\subsection{1-point function on the torus}\label{s2.3}
In order to motivate Theorem 2, consider the $T\overline T$ deformation of the one-point function of a local operator on the torus, which has a natural action of $\Gamma_1=$SL$(2,{\mathbb Z})$. A 2-torus may be thought of as ${\mathbb R}^2/({\mathbb Z}R_a+{\mathbb Z}R_b)$, where $R_a,R_b\in{\mathbb R}^2$, such that the area $R_a\wedge R_b>0$. The generators of $\Gamma_1$ are $S:(R_a,R_b)\to(R_b,-R_a)$ and 
$T:(R_a,R_b)\to(R_a,R_b+R_a)$. Rotating to a basis where $R_a=(|R_a|,0)$, we can write $R_b=|R_a|(-\d_2,\d_1)$, where the usual modulus is $\tau=\i\d$ where $\d=\d_1+\i\d_2$ with $\d_1>0$.

In a translationally invariant field theory, the 1-point function $\langle\Phi(x)\rangle$ of a scalar operator $\Phi(x)$ on the torus is independent of $x$, and in a CFT has the form
$|R_a|^{-k}F_2^0(\d)$, where now $k$ is the scaling dimension of $\Phi$. The symmetry under $\Gamma_1$ implies that $F_2^0(1/\d)=|\d|^kF_2^0(\d)$ and $F_2^0(\d+\i)=F_2^0(\d)$. Moreover, in a CFT,   $F_2^0(\d)$ has a Fourier expansion as in (\ref{1.7}).\footnote{More generally, a sum over several $\Delta$s.} 

Paralleling the discussion in the previous section, we define the double Laplace transform
\be
\Omega^0(R_a,s)=\int_{R_a\wedge R'_b>0}\erm^{-s.R_b'}\langle\Phi\rangle^0(R_a,R'_b)d^2R'_b\,.
\ee
\be
\langle\Phi\rangle^0(R_a,R_b)=\int \erm^{s.R_b}\Omega^0(R_a,s)(d^2s/(2\pi\i)^2)\,.
\ee

Choosing $x=X$, and running the contour $C$ in (\ref{CXx}) from $X$ to $X$ around an $a$ cycle, the effect of the $T\overline T$ deformation is to send $R_a\to R_a-\l\wedge N_a$, where $N_a\sim -s$ is the force acting across $C$.  Thus
\be
\langle\Phi\rangle^\l(R_a,R_b)=\int \erm^{s.R_b}\Omega^0(R_a-\l\wedge s,s)(d^2s/(2\pi\i)^2)\,.
\ee
We note in passing that, at least formally, this implies the PDE
\be\label{pde2}
\p_\l\langle\Phi\rangle^\l(R_a,R_b)=-(\p_{R_a}\wedge\p_{R_b})\langle\Phi\rangle^\l(R_a,R_b)\\,
\ee
which is $\Gamma_1$ invariant. 

In a CFT,
\be
\Omega^0(R_a,s)=\int \erm^{-s.{\bf d}'.R_a}|R_a|^{-k+2}F^0({\bf d}')d^2{\bf d}'\\,
\ee
where $R'_b={\bf d'}.R_a$, i.e.
\be
\left(\begin{array}{c}R_b^1\\R_b^2\end{array}\right)=
\left(\begin{array}{cc}\d_2 & -\d_1 \\ \d_1 & \d_2 \end{array}\right)
\left(\begin{array}{c}R_a^1\\R_a^2\end{array}\right)
\ee
and similarly $R_b={\bf d}.R_a$.
Then
\be
\langle\Phi\rangle^\l(R_a,{\bf d}.R_a)=\int \erm^{s.{\bf d}.R_a}\int \erm^{-s.{\bf d}'.(R_a-\l\wedge s)}
|R_a-\l\wedge s|^{-k+2}F({\bf d}'))d^2{\bf d}'(d^2s/(2\pi\i)^2)\,.
\ee
In terms of components, in a frame where $R_a^2=0$, this is
\be
\int\int[(|R_a|-\l s_1)^2+\l^2s_2^2]^{-k/2+1}
\erm^{\l\d'_1(s_1^2+s_2^2)+s_1|R_a|(\d_1-\d'_1)+s_2|R_a|(\d_2-\d'_2)}
F^0(\d')d^2d'(d^2s/(2\pi\i)^2)\\,
\ee
which should be compared to (\ref{eq}). Defining $\a=\l/$(area) $=\l/(\d_1|R_a|^2)$, and rescaling $s$, we find
\be
F_2^\a(\d)=\int_C\int_{\mathbb H}[(1-\a\d_1 s_1)^2+\a^2\d_1^2s_2^2]^{-k/2+1}
\erm^{\a\d_1\d'_1s^2+s.(\d-\d')}
F_2^0(\d')d^2\d'(d^2s/(2\pi\i)^2)\,,
\ee
where the  $s_{1,2}$ contours $C$ lie up the imaginary axis, and $\d'$ is integrated over the right half plane.
As before, if $\Delta>0$  this converges uniformly for $\d$ in any closed subset of $\mathbb H$, and although the arguments leading up to this lack rigor, we may now take it as a definition of the deformed $F_2$, and thence the deformed
one-point functions.

\subsubsection{Partition function}\label{s2.3.1}
The CFT torus partition function $Z(R_a,R_b)$ has the form (\ref{1.7}), but is modular invariant. However, its deformed version is not given by (\ref{e2}) with $k=0$, because the deformation should be applied to the partition function with a marked point $X$ where $u(X)=0$,  that is $Z$ divided by the area $R_a\wedge R_b=\d_1|R_a|^2$. 

As a result the above arguments are slightly modified: (\ref{eq}) becomes
\be
{\d_1}^{-1}Z^\l(R_a,{\bf d}.R_a)=\int \erm^{s.{\bf d}.R_a}\int \erm^{-s.{\bf d}'.(R_a-\l\wedge s)}
{\d'_1}^{-1}F({\bf d}'))d^2{\bf d}'(d^2s/(2\pi\i)^2)\,,
\ee
so that
\be
Z^\a(\d)=\int\int
\erm^{\a\d_1\d'_1s^2+s.(\d-\d')}
(\d_1/\d'_1)Z^0(\d')d^2\d'(d^2s/(2\pi\i)^2)\,
\ee
which again may be used as a definition of the deformed partition function.

\section{Proof of Theorem 1}\label{s3}
Motivated by the above discussion, given a function $F^0(\d)=\sum_{n=0}^\infty a_nq^{\Delta+n}$ where $q=\erm^{-2\pi\d}$ with Re\,$\d>0$, which satisfies $F^0(1/\d)=\d^kF^0(\d)$, we define its $T\overline T$ deformation by the limit as $\e\to0$, if it exists, of
\be\label{def2}
F_\e^\a(\d)=\int_{-\i\infty}^{\i\infty}\erm^{s\d}\int_\e^{1/\e} (1-\a\d s)^{1-k}\erm^{-s\d'(1-\a\d s)}F^0(\d')d\d'(ds/2\pi\i)\\,
\ee
where the branch cut of $z^{1-k}$ is taken to lie along the negative real $z$ axis. As long as $\a>0$ and Re\,$\d>0$, the integral is uniformly convergent and it is permissible to interchange the orders of integration.

Inserting the uniformly convergent expansion of $F^0(\d')$, we may interchange the orders of summation and integration to find 
\be
F_\e^\a(\d)=\sum_{n=0}^\infty a_n\int_{-\i\infty}^{\i\infty}\int_\e^{1/\e} (1-\a\d s)^{1-k}
\erm^{s(\d-\d')+\a\d\d' s^2-2\pi(\Delta+n)\d'}d\d'(ds/2\pi\i)\,.
\ee
As long as $\Delta>0$ we can take the limit $\e\to0$ term by term, giving
\be\label{3.2}
\sum_{n=0}^\infty a_n\int_{-\i\infty}^{\i\infty} \frac{\erm^{s\d}(1-\a\d s)^{1-k}}
{2\pi(\Delta+n)+s-\a\d s^2}\frac{ds}{2\pi\i}\,.
\ee
The integrand has simple poles at $s=s_\pm=(1/2\a\d)(1\pm\sqrt{1+8\pi(\Delta+n)\a\d})$ and a branch cut from $s=(\a\d)^{-1}$ to $+\infty$. For $\Delta>0$ only $s_-$ lies to the left of the $s$ integration contour, and we may pull it back and evaluate the residue to obtain the expression in (\ref{e1}).

On the other hand, we may first perform the $s$ integral by completing the square in the exponent, writing it as $\a\d\d'(s-(\d'-\d)/2\d\d')^2-(\d'-\d)^2/4\a\d\d'$.
Setting $s=(\d'-\d)/2\d\d'+\i (\d\d')^{-1/2}t$ and shifting the contour so that $t$ is real,
we find, after some algebra,
\be\label{fa}
F_\e^\a(\d)=\int_{-\e}^\e K^\a(\d,\d')(\d/\d')^{k/2}F^0(\d')(d\d'/\d')\,,
\ee
where
\be
K^\a(\d,\d')=\erm^{-(\d'-\d)^2/4\a\d\d'}\int_{-\infty}^\infty\big((\d+\d')/2(\d\d')^{1/2}-\i t\big)^{1-k}\erm^{-\a t^2}dt\,.
\ee
If $\Delta>0$ we may remove the $\e$ cutoff in (\ref{fa}). The theorem now follows on recognizing that both $K^\a(\d,\d')$ and the measure $d\d'/\d'$ are invariant under $(\d,\d')\to(1/\d,1/\d')$, which implies that if $\d^{-k/2}F_1^0(\d)$ is invariant so also is ${\d}^{-k/2}F_1^\a(\d)$. $\square$

\subsubsection*{Remarks.}
\begin{itemize}
\item If $\Delta<0$ then (\ref{e1}) converges for real $\d$ as long as $(4\a\d)^{-1}>2\pi|\Delta|$ and
$(4\a/\d)^{-1}>2\pi|\Delta|$, that is $(8\pi\Delta\a)^{-1}<\d<8\pi\Delta\a$. As long as this holds, $s_-$ and $s_+$ are both real and we may shift the contour in (\ref{3.2}) so as to lie between them. Thus the theorem still holds in this restricted domain. $F^\a(\d)$ has square root singularities at the end points, which in the physics literature are called Hagedorn singularities. 
\item the proof of $S$-invariance depends only on the invariance of the kernel $K^\a(\d,\d')$. Thus we could use some other such kernel. However this will not in general lead to a point power spectrum in (\ref{e1}) for $\a>0$.
\item If $\a<0$ we can still define $F^\a$ by wrapping the contour in (\ref{def2}) to lie just above and below the real axis so as to include all the poles at $s=-2\pi(\Delta+n)$ when $\a=0$.
Then for $\a<0$ it includes the finite number of  poles $s_\pm$ on the real axis, as well as the branch cut which now runs from $-\infty$ to  $-(|\a|\d)^{-1}$, but none of the infinity of poles with a non-zero imaginary part. Thus the spectrum of powers of $q$ has a continuum  as well as a finite discrete part. The proof of $S$-invariance of ${\d}^{-k/2}F_1^\a(\d)$ then proceeds as before, picking up only the real spectrum. However the partition function has a square root singularity every time a pair of roots $s_\pm$ meet and become complex. 
\end{itemize}

\section{Deformed Mellin transform and Dirichlet series}\label{s4}
To any function of the form $F(\d)=\sum_{n=0}^\infty a_nq^{\Delta+n}=\sum_na_n\erm^{-2\pi(\Delta+n)\d}$ with $\Delta>0$
there corresponds a Dirichlet series
\be
\phi(s)=\sum_{n=0}^\infty\frac{a_n}{(\Delta+n)^s}=
\frac{(2\pi)^s}{\Gamma(s)}R(s)\,,
\ee
where 
$R(s)=\int_0^\infty\d^{s-1}F(\d)d\d$ is the Mellin transform of  $F$.
Moreover if $F(\d)=\d^{-k}F(1/\d)$, then $F(\d)$ grows no worse than $\d^{-k}$ as $\d\to0$, so $R(s)$ is analytic for  Re\,$s>k$, and admits an analytic continuation satisfying
\be
R(s)=R(k-s)\,.
\ee 
\gdef\b{{\beta}}
However, as may be seen from (\ref{e1}), the deformed function $F^\a(\d)$, at fixed $\a$, is not a sum of powers of $q=e^{-2\pi\d}$, although it is at fixed $\b=\a\d$. Thus there are two alternative candidates for the deformed $\phi(s)$.

\subsection{Deformed  Mellin transform}
The simplest and most interesting definition is through the Mellin transform
$\phi^\a(s)=\frac{(2\pi)^s}{\Gamma(s)}R^\a(s)$ where $R^\a(s)=\int_0^\infty\d^{s-1}F^\a(\d)d\d$. 

We now outline the elements of the proof of Thm.~1a. By definition,
\be
R^\a(s)=\int_0^\infty\d^{s-1}\int_{-\i\infty}^{\i\infty}\erm^{t'\d}\int (1-\a\d t')^{1-k}\erm^{-t'\d'(1-\a\d t')}F^0(\d')d\d'(dt'/2\pi\i)d\d\,.
\ee
This looks complicated, but on substituting $\d=u\d'$  and rescaling
$t'\to t'/\d'$,
\begin{eqnarray}
R^\a(s)&=&\int_0^\infty(u\d')^{s-1}\int_{-\i\infty}^{\i\infty}\erm^{t'u}\int (1-\a ut')^{1-k}\erm^{-t'(1-\a ut')}F^0(\d')d\d'(dt'/2\pi\i) du\\
&=&I^\a(k;s)\,R^0(s)\,,
\end{eqnarray}
where
\be\label{Idef}
I^\a(k;s)=\int_0^\infty u^{s-1}\int_{-\i\infty}^{\i\infty} (1-\a ut)^{1-k}\erm^{tu-t(1-\a ut)}(dt/2\pi\i)\, du
\ee
Now we follow the same method as in the proof of Thm.~1, completing the square in the exponent, to find, after some algebra,
\be
I^\a(k;s)=\int_0^\infty u^{s-k/2-1}\erm^{-(u-1)^2/4\a u} \int_{-\infty}^{\infty}(\i\a t+(u^{1/2}+u^{-1/2})/2)^{1-k}\erm^{-\a t^2}(dt/2\pi) du\,.
\ee
The integrals are absolutely uniformly convergent and define an entire function of $s$.
Under $u\to1/u$, $u^{s-k/2-1}du\to u^{-s+k/2-1}du=u^{k-s-k/2-1}du$, so $I^\a(k;s)=I^\a(k;k-s)$.     $\square$

In fact $I^\a(k;s)$ is proportional to a confluent hypergeometric function.
Going back to (\ref{Idef}) and rescaling $t\to t/\a u$, then $u\to1/u$,
\be
I^\a(k;s)=\a^{-1}\int_0^\infty u^{-s}\int_{-\i\infty}^{\i\infty} (1-t)^{1-k}\erm^{(t/\a)-u(t(1-t)/\a)}(dt/2\pi\i)\,du
\ee
\be
=\a^{-s}\Gamma(1-s)\int_{-\i\infty}^{\i\infty} \erm^{t/\a}t^{s-1}(1-t)^{s-k}(dt/2\pi\i)
\ee
\be
=\a^{-s}\Gamma(1-s)\pi^{-1}\sin\pi(s-k)\int_0^1 \erm^{t/\a}t^{s-1}(1-t)^{s-k}dt
\ee
\be
=\a^{-s}\frac{\sin\pi(s-k)\Gamma(s-k+1)}{\sin\pi s\Gamma(2s-k+1)}{}_1F_1(s,2s-k+1;1/\a)\,.
\ee

\subsection{Deformed Dirichlet series}
Alternatively, we may define the $T\overline T$-deformed Dirichlet series by 
\be
\phi^\b(s)=\sum_{n=0}^\infty\frac{a_n^\b}{\big((\Delta+n)^\b\big)^s}\,,
\ee
where 
 $a_n^\b$ and $(\Delta+n)^\b$ are the deformed coefficients and exponents defined in (\ref{e1},\ref{xb}), with $\a=\b/\d$.
 As before, we have $\phi^\b(s)=\frac{(2\pi)^s}{\Gamma(s)}R^\b(s)$ where $R^\b(s)=\int_0^\infty\d^{s-1}F^\b(\d)d\d$, the only difference being that the integration is performed at fixed $\b$ rather than fixed $\a$. Thus
\be
R^\b(s)=\int_0^\infty\d^{s-1}\int_{-\i\infty}^{\i\infty}\erm^{t'\d}\int (1-\b t')^{1-k}\erm^{-t'\d'(1-\b t')}F^0(\d')d\d'(dt'/2\pi\i)d\d\,.
\ee
The $\d$ integral is now immediate, and $F^0$ may be expressed in terms of $R^0$ by an inverse Mellin transform. 
\begin{eqnarray}
R^\b(s)&=&\Gamma(s)\int_{-\i\infty}^{\i\infty}\int_0^\infty\int_{-\i\infty}^{\i\infty}{t'}^{-s}(1-\b t')^{1-k}\erm^{-t'\d'(1-\b t')} {\d'}^{-s'}R^0(s')\frac{ds'}{2\pi\i}d\d'\frac{dt'}{2\pi\i}\\
&=&\Gamma(s)\int_{-\i\infty}^{\i\infty}\int_{-\i\infty}^{\i\infty}{t'}^{s'-s-1}(1-\b t')^{s'-k}\Gamma(1-s')R^0(s')\frac{ds'}{2\pi\i}\frac{dt'}{2\pi\i}\,.
\end{eqnarray}
This has an asymptotic expansion in powers of $\b$ of the form
\be
R^\b(s)=\sum_{r=0}^\infty c_r\b^rR^0(s-r)\,,
\ee
which shows that reflection symmetry under $s\to k-s$ is lost once $\beta\not=0$, as expected, since the $\d\to1/\d$ symmetry of $F$ holds only at fixed $\a$.

\section{Proof of Theorem 2}\label{s5}
Motivated now by the discussion in Sec.~(\ref{s2.3}), given a real function $F_2^0(\d)=$\\ $\sum_{n=0}^\infty \sum_{\bar n=0}^\infty a_{n,\bar n}q^{\Delta+n}{\bar q}^{\Delta+\bar n}$ where $q=\erm^{-2\pi\d}$ with Re\,$\d>0$, with $a_{n,\bar n}=a_{\bar n,n}$, which satisfies $F_2^0(1/\d)=|\d|^kF_2^0(\d)$ and
$F_2^0(\d+\i)=F_2^0(\d)$, we define its $T\overline T$ deformation by 
\be\label{f2def}
F_2^\a(\d)=\int\int[(1-\a\d_1 s_1)^2+\a^2\d_1^2s_2^2]^{-k/2+1}
\erm^{\a\d_1\d'_1\vec s^2+\vec s.(\vec\d-\vec\d')}
F_2^0(\d')d^2\d'(d^2s/(2\pi\i)^2)\,,
\ee
where we again assume that $\Delta>0$ so the integrals converge uniformly.
Substituting the convergent expansion (\ref{1.7}), in each term $\propto a_{n,p}$ integrating on $\d'_2$ sets $s_2=2\pi p$, and then the $\d'_1$ integral gives
\be\label{5.2}
\int\frac{[(1-\a\d_1 s_1)^2+4\pi^2\a^2\d_1^2p^2]^{-k/2+1}\erm^{s_1\d_1}}
{2\pi(\Delta+n)+4\pi^2\a^2\d_1^2p^2+s_1-\a\d_1s_1^2}\frac{ds_1}{2\pi\i}\,,
\ee
now with simple poles at $s_1=s_\pm=(1/2\a\d_1)(1\pm\sqrt{1+8\pi(\Delta+n)\a\d_1+16\pi^2\a^2\d_1^2p^2})$. Again evaluating the residue at $s_-$ then gives (\ref{e2}).

On the other hand, completing the square in the exponent and setting
$\vec t=\vec s+(\vec\d-\vec\d')/2\a\d_1\d'_1$
we find, after some algebra,
\be
F_2^\a(\d)=\int_{\mathbb H}K_2^\a(\d,\d')(\d'_1/\d_1)^{k/2}F_2^0(\d)(d^2\d'/{\d'_1}^2)\,,
\ee
where
\begin{eqnarray}
K_2^\a(\d,\d')&=&\a^{-1}\erm^{-(\vec\d-\vec\d')^2/4\a\d_1\d'_1}\times\nonumber\\
&&\int\left[\left(\frac{(\d_1+\d'_1)}{2(\d_1\d'_1)^{1/2}}-t_2\right)^2+\left(\frac{(\d_2-\d'_2)}{2(\d_1\d'_1)^{1/2}}- t_1\right)^2\right]^{-k/2+1}\frac{\erm^{-\vec t^2}d^2t}{(2\pi\i)^2}\,.
\end{eqnarray}
Note that the exponent in the first line is the square of the hyperbolic distance between $\d$ and $\d'$.
At this point we can invoke the reality of $F^0_2$ to change the sign of $\d'_2$ in the second term. The integrand and the measure are then invariant under SO$(2)$ rotations, so we may rotate to a frame where $\d_2+\d'_2=0$. Thus the integral is equal to
\be
\int\left[\left(\frac{|\d+\d'|}{2(\d_1\d'_1)^{1/2}}-t_2\right)^2+ t_1^2\right]^{-k/2+1}\frac{\erm^{-\vec t^2}d^2t}{(2\pi\i)^2}\,.
\ee
\gdef\db{{\bar\delta}}
We then have, as before, $K_2^\a(1/\d,1/\d')=K_2^\a(\d,\d')$. This is most easily seen in terms of $\d=\d_1+\i\d_2$ and $\db=\d_1-\i\d_2$, since 
\be
 \frac{|\vec\d\pm\vec\d'|^2}{\d_1\d'_1}=\frac{(\d\pm\d')(\db\pm\db')}{\d_1\d'_1}\to
 \frac{((\d'\pm\d)/\d\d')(\db'\pm\db)/\db\db')}{(\d_1/\d\db)(\d'_1/\d'\db')}
 =\frac{(\d'\pm\d)(\db'\pm\db)}{\d_1\d'_1}\,.
 \ee
Since the measure $d^2\d'/{\d'_1}^2$ is invariant we conclude that if $\d_1^{k/2}F^0(\d)$ is invariant under $\d\to1/\d$, then so is $\d_1^{k/2}F^\a(\d)$. This establishes Theorem 2, since $\d_1\to\d_1/|\d|^2$. $\square$

Again, the above argument is  strictly valid only if $\Delta>0$: otherwise  we should restrict the range of $\d$ as before.

\subsection{Proof of Theorem 2a}\label{s5.1}
As was discussed in Sec.~\ref{s2.3.1}, if $k=0$, that is $F_2^0$ is $\Gamma_1$ invariant, there is an alternate version in which
\be\label{e34}
F_2^\a(\d)_\e=\int\int_\e \erm^{\a\d_1\d'_1s^2+s.(\d-\d')}(\d_1/\d'_1)
F_2^0(\d')d^2\d'(d^2s/(2\pi\i)^2)\,,
\ee
where $\int_\e$ indicates that the $\d'_1$ integration is over $(\e,\e^{-1})$. 

As compared with (\ref{5.2}), after substituting the expansion (\ref{1.7}), the $\d'_1$ integration gives
\be
\log[(2\pi(\Delta+n)+4\pi^2\a^2\d_1^2p^2+s_1-\a\d_1s_1^2)/\e]+O(1)
=\log(s_1-s_-)+\log(s_+-s_1)+O(1)\,,
\ee
where the remainder is either non-singular in $s_1$ or vanishes as $\e\to0$.
There are now branch cuts along $(-\infty,s_-)$ and $(s_+,+\infty)$ and the contour in $s_1$ separates them. Wrapping this around the left hand cut then gives
$\int_{-\infty}^{s_-}\d_1\erm^{s_1\d_1}ds_1=\erm^{s_-\d_1}$, giving (\ref{e2a}) with no prefactor.

On the other hand, the $(s_1,s_2)$ integrations can now be performed explicitly, to give
\be
F_2^\a(\d)_\e=(4\pi\a)^{-1}\int_\e \erm^{-|\d-\d'|^2/4\a\d_1\d'_1}
F_2^0(\d')(d^2\d'/{\d'_1}^2)\,,
\ee
which is equivalent to the result for the deformed CFT partition function first obtained by Dubovsky, Gorbenko and Hern\'andez-Chifflet \cite{DubJT}. Since the kernel and the integration measure are $\Gamma_1$ invariant, this establishes Theorem 2a. $\square$

\subsubsection*{Remarks.}
\begin{itemize}
\item The theorem extends to a multiplet of functions $\{F^\a_I\}$ transforming linearly according to some representation of  $\Gamma_1$, that is 
\be
F^\a_I(1/\d)=\sum_J{\bf\rm S}^I_JF^\a_I(\d)\,.
\ee
If the $\{F^0_I\}$ satisfy this, so also do the $\{F^\a_I\}$. An example is given the characters
$\chi_i(q)\chi_{\bar i}(\bar q)+\chi_i(\bar q)\chi_{\bar i}(q)$ of the product ${\cal V}\otimes\overline{\cal V}$ of two Virasoro algebras arising in the decomposition of the torus partition function, recently analyzed in detail in \cite{dattach}.
\item It is tempting to try the same variant construction for the holomorphic case of Theorem 1, for example the Virasoro characters $\{\chi_i\}$ themselves, maintaining the coefficients $a_n$ with no prefactor. However, the factors of $\d'$ do not arrange themselves as conveniently as in (\ref{e34}), and as a consequence, if $\chi^0_i(1/\d)=\sum_jS^i_j\chi^0_j(\d)$, this no longer holds for the deformed versions. [This would correspond to a deformed CFT partition function on an annulus, which has no reason to be invariant.]
\item  The situation when $\a<0$ is more tricky than in Thm.~1. There is an infinite number of singularities $s_\pm$ with $n\sim|p|\gg1$ which remain almost undeformed on the real axis, and any choice of the contour which includes the undeformed singularities will necessarily include an infinity  of these. However, with such a contour, the above argument for modular invariance of the deformed partition function at fixed $\a$ still goes through, although it now has a dense set of singularities as $\a$ is varied and the singularities at $s_\pm$ pinch and go into the complex plane. 
\end{itemize}

\section{Deformed Maass forms}\label{s6}
Recall that a \em Maass form \em for $\Gamma_1 =$ SL$(2,{\mathbb Z})$ is a smooth function $F_M$ on $\mathbb H$ satisfying
the following three conditions:
\begin{enumerate}
\item[(i)] for all $\gamma\in\Gamma_1$, $F_M(\gamma(\d))=F_M(\d)$;
\item[(ii)] $F_M$ is an eigenfunction of the invariant Laplacian on the fundamental region $\Delta_{\cal F}=-\d_1^2(\p_{\d_1}^2+\p_{\d_2}^2)$;
\item[(iii)] $F_M(\d_1+\i\d_2)$ is polynomially bounded as $\d_1\to\infty$.
\item[(iv)] A Maass form is a Maass \em cusp form \em if also 
 its Fourier coefficient $f_0$ in $F_M=\sum_{p\in{\mathbb Z}}f_p(\d_1)\erm^{2\pi\i p\d_2}$ vanishes.
 \end{enumerate}
 That Maass forms are defined in this somewhat abstract manner rather than by a $q$-series makes it less obvious how to define their $T\overline T$ deformation. However, we have the following 
 
 \noindent{\bf Lemma.} \em If $F(\d=-\i R_b/R_a)$ is a smooth function on the torus which depends only on the modulus $\d$ (not necessarily holomorphically), then\em
 \be
 A(\p_{R_a}\wedge\p_{R_b})\,F=(1/4)\Delta_{\mathbb H}\,F\,.
 \ee
 The proof is straightforward, using $A=(\i/2)(R_aR_b^*-R_a^*R_b)$ and
 $\p_{R_a}\wedge\p_{R_b}=-(\i/2)(\p_{R_a}\p_{R_b^*}-\p_{R_a^*}\p_{R_b})$.
 In Sec.~\ref{s2.3} we argued that the $T\overline T$ variation of any such function in a CFT is given by 
 \be
 \p_\l F|=-(\p_{R_a}\wedge\p_{R_b})F=-(1/4A)\Delta_{\mathbb H}\,F\,,
 \ee
 so that, writing $\l=A\a$, $\p_\a F=-\frac14\Delta_{\mathbb H}\,F$. Thus any eigenfunction of the Laplacian which depends only on the modulus, is multiplicatively transported by the $T\overline T$ flow. If $F$ is also $\Gamma_1$ invariant, we may restrict to the fundamental domain $\cal F$. Thus, for a Maass form of eigenvalue $\Lambda$,
 \be
  F_M^{\a}=\erm^{-(\Lambda/4)\a}\,F^0_M\,.
  \ee
  This completes the proof of Thm.~3. $\square$
  
  The $L^2({\cal F})$ eigenfunctions and eigenvalues of $\Delta_{\cal F}$  are well characterized \cite{ter}. Aside from the constant function a complete set is given by a continuous spectrum, the real Eisenstein series $E_s$ (see below) with Re\,$s=\frac12$ and $\Lambda=s(1-s)$, which have power law decay as $\d_1\to\infty$, and a discrete series of somewhat intractable Maass cusp forms which decay exponentially with a large gap. Thus the $T\overline T$ deformation of a generic function $F\in L^2({\cal F})$,
  for example $F_2^\a$ in (\ref{e2}) with $\Delta>0$, should decay to a constant value exponentially at a rate $\propto\erm^{-\a/16}$.

 \subsection{Deformed non-holomorphic Eisenstein series}\label{s6.2}
Important examples of Maass forms are the non-holomorphic Eisenstein series
\be
E_s(\d)=\sum_{(m,n)\in{\mathbb Z}^2\setminus(0,0)}\frac{\d_1^{s}}{|\i m\d+n|^{2s}}\,,
\ee
which may be thought of as a uniform sum over $\Gamma_1$ of images of the $(0,1)$ term, so is automatically invariant under the group action. Also, since $\Delta_{\mathbb H}$ is invariant, its action on each term in the sum is similar. On the $(0,1)$ term we have simply $-\d_1^2\p_{\d_1}^2\d_1^s=s(1-s)\d_1^s$, so $\Lambda=s(1-s)$.

On the other hand, we can think of $E_s$ as a form on a dimensionful torus as
\be
E_s=\sum_{(m,n)\in{\mathbb Z}^2\setminus(0,0)}\frac{(A/|R_a|^2)^{s}}{|\i m(R_b/R_a)+n|^{2s}}
=\sum_{(m,n)\in{\mathbb Z}^2\setminus(0,0)}\frac{A^{s}}{|\i mR_b+nR_a|^{2s}}\,.
\ee
By rotational symmetry and scaling, the $(m,n)$ term deforms into some function  
\be
f_s^\a\left(\frac{A}{|\i mR_b+nR_a|^2}\right)=f_s^\a\left(\frac{R_a\wedge R_b}{|\i mR_b+nR_a|^2}\right)\,,
\ee
satisfying
\be
\p_\a f_s^\a=-A(\p_{R_a}\wedge\p_{R_b})f_s^\a=-(1/4)\Delta_{\mathbb H}f_s^\a\,,
\ee
using the Lemma.
Again choosing the $(0,1)$ term as representative, this becomes simply
\be
\p_\a f_s^\a(X)=-(1/4)X^2{f_s^\a\,}''(X)\,,
\ee
with the initial condition $f^0_s(X)=X^s$, and the solution 
\be
f^\a_s(X)\propto e^{-(1/4)s(1-s)\a}X^s\,,
\ee
in agreement with Thm.~3.

\section{Some examples}\label{s7}
\subsection{Deformed Jacobi theta functions}
A simple example of Thm.~1 is that
the well known identity
\be
\vartheta_3(0;\d)=\sum_{n\in{\mathbb Z}}\erm^{-\pi n^2\d}=\d^{-1/2}\sum_{n\in{\mathbb Z}}\erm^{-\pi n^2/\d}=\d^{-1/2}\vartheta_3(0;1/\d)
\ee
becomes a similar identity for
\be\label{thetadef}
\vartheta^\a_3(0;\d)\equiv
\sum_{n\in{\mathbb Z}}\frac{\sqrt{(1+\sqrt{1+4\pi\a n^2\d})/2}}{\sqrt{1+4\pi\a n^2\d}}
\erm^{-(1/2\a)(\sqrt{1+4\pi\a n^2\d}-1)}\,,
\ee
which can readily be checked numerically.

More generally, the symmetrized version of the inversion relation for the Jacobi theta function
\be
\vartheta_3(z\d^{1/2};\d)\equiv\sum_n\erm^{-\pi n^2\d+2\pi\i nz\d^{1/2}}
=\d^{-1/2}\erm^{-\pi z^2}\vartheta_3(\i z/\d^{1/2};1/\d)
\ee
is also satisfied by its deformed version, given by substituting 
$\sqrt{1+4\pi\a n^2\d}\to$\\$\sqrt{1+4\pi\a n^2\d-8\pi\i nz\d^{1/2}}$ in (\ref{thetadef}). 

Similarly, a Jacobi form of weight $k$ and index $m$ obeys 
\be
\phi(z\d^{1/2};\d)
=\d^{-k}\erm^{-\pi mz^2}\phi(\i z/\d^{1/2};1/\d)\,,
\ee
as does its deformed version, defined by further modifying the power in the numerator of the prefactor (\ref{thetadef}) in from $\frac12$ to $(1-k)/2$.

\subsection{Deformed partition sums}
Let $P(n)$ be the number of distinct partitions of $n\in{\mathbb N}^+$ into positive integers.
Since $\eta(\tau)^{-1}=q^{-1/24}\sum_{n=0}^\infty P(n)q^n$ is a form with $k=-\frac12$, by Thm.~1
\begin{eqnarray}
&&\sum_{n=0}^\infty P(n)\frac{(1+\sqrt{1+8\pi\a(n-1/24)\d})^{3/2}}{\sqrt{1+8\pi\a(n-1/24)\d}}
\erm^{-(1/2\a)\sqrt{1+8\pi\a(n-1/24)\d}}\nonumber\\
&&=\d^{1/2}
\sum_{n=0}^\infty P(n)\frac{(1+\sqrt{1+8\pi\a(n-1/24)/\d})^{3/2}}{\sqrt{1+8\pi\a(n-1/24)/\d}}
\erm^{-(1/2\a)\sqrt{1+8\pi\a(n-1/24)/\d}}\,.
\end{eqnarray}
In this case $\Delta<0$, and so this is valid only for $\pi\a/3<\d<3/\pi\a$. The $n=0$ term on the right hand side is singular at $\d/\a=\pi/3$, dictating the radius of convergence of the left hand side, and \em vice versa\em. The singularities of the right hand side are of square root type rather than essential as in the undeformed case, and should determine the Hardy-Ramanujan asymptotics of $P(n)$.  It would be interesting to understand how these contrive to be independent of $\a$.

\subsection{Deformed Eisenstein series}
The undeformed series
\be
E_k(\d)=\sum_{(m,n)\not=(0,0)}\frac1{(m+\i n\d)^k}
\ee
gives an example of a modular form of weight $k$ (strictly only if $k$ is a positive even integer).
To use (\ref{def2}) is difficult. But if we consider 
\be
Z(L_1,L_2)\equiv L_1^{-k}E_k(L_2/L_1)=\sum_{(m,n)\not=(0,0)}\frac1{(mL_1+\i nL_2)^k}\,,
\ee
it satisfies $Z(L_1,L_2)=Z(L_2,L_1)$ and it can be shown that the scaling solution of the PDE
\be
\p_\l Z^\l(L_1,L_2)=-\p_{L_1}\p_{L_2}Z^\l(L_1,L_2)
\ee
is $E_k^\a(\d)$ with $\a=\l/(L_1L_2)$ and $\d=L_2/L_1$. For each $(m,n)$ we then look for a solution of the form $f^\l(L=mL_1+\i nL_2)$, so that
$\p_\l f=-\i mn\p_L^2f$, which can be solved by Green function. The final result is
\be
E_k^\a(\d)=(4\pi\a)^{-1/2}\sum_{(m,n)\not=(0,0)}\int_{-\infty}^\infty \frac{\erm^{-\ell_{m,n}^2/4\a}}{(\ell_{m,n}(-\i  mn\d)^{1/2}+m+\i  n\d)^k}d\ell_{m,n}\,,
\ee
so that each lattice point gets an independent gaussian deformation, which vanishes on the axes $m=0$ and $n=0$, and on the diagonals $m=\pm n$ is purely transverse. This preserves the symmetry under $S$ but not $T$, as expected. (We showed in Sec.~\ref{s6.2} that the non-holomorphic version, which gives a Maass form, does not evolve, up to a multiplicative constant.)

\subsubsection*{Acknowledgements.} This work was begun under remote participation in the KITP program on Modularity in Quantum Systems in Fall 2020, supported in part by the National Science Foundation under Grant No. NSF PHY-1748958, and continued  with support from the Quantum Science Center (QSC), at the University of California, Berkeley, a National Quantum Information Science Research Center of the U.S. Department of Energy (DOE).

\end{document}